\input amstex
\documentstyle{amsppt}
\input bull-ppt
\keyedby{bull293e/PAZ}

\define\lan{\langle}
\define\ran{\rangle}

\topmatter
\cvol{27}
\cvolyear{1992}
\cmonth{July}
\cyear{1992}
\cvolno{1}
\cpgs{139-142}
\title A class of nonsymmetric harmonic Riemannian
spaces \endtitle
\author Ewa Damek and Fulvio Ricci \endauthor
\address Instytut Matematyczny Universytetu Wroc\l awskiego,
 pl. Grunwaldzki 2/4, 50-384 Wroc\l aw, Poland\endaddress
\address Dipartimento di Matematica, Politecnico di 
Torino, Corso
Duca degli Abruzzi 24, 10129 Torino, Italy\endaddress
\ml FRICCI\@POLITO.IT\endml
\subjclass Primary 53C25, 53C30, 43A85, 22E25, 
22E30\endsubjclass
\date July 11, 1991\enddate
\abstract Certain solvable extensions of $H$-type
groups provide noncompact counterexamples to a
conjecture of Lichnerowicz, which asserted that
``harmonic'' Riemannian spaces must be rank 1 symmetric
spaces.\endabstract
\endtopmatter

\document

A Riemannian space $M$ with Laplace-Beltrami operator
$\Delta$ is called {\it harmonic} if, given any function
$f(x)$ on $M$ depending only on the distance $d(x,x_0)$
from a given point $x_0$, then also $\Delta f(x)$ depends
only on  $d(x,x_0)$. 

Equivalently, $M$ is harmonic if for every $p\in M$ the
density function $\omega_{x_0}(x)$ expressed in terms of
the normal coordinates around the point $x_0$ is a
function of $d(x,x_0)$ (see~\cite{1,~11}). 

In 1944  Lichnerowicz~\cite{10} proved that in dimensions
not greater than 4 the harmonic spaces coincide with the
rank-one symmetric spaces. He also raised the question of
determining whether the same is true in higher dimensions.

Among the most recent progress made on the so-called
Lichnerowicz conjecture, Szab\'o~\cite{11} proved it to hold
true in arbitrary dimension for compact manifolds with 
finite fundamental group.

In this announcement we present a counterexample that
arises in the noncompact case. It proves the
Lichnerowicz conjecture not to be true in general for
infinitely many dimensions, the smallest of them being 7.

This example is based on the notion of H-type group due
to Kaplan~\cite{8}, and on the geometric properties of
their one-dimensional extensions $S$ introduced by
Damek~\cite{5} and studied also in~\cite{2,~3,~4,~6}.

The rank-one symmetric spaces of the noncompact type
different from the hyperbolic spaces are special examples
of such groups $S$. Even though symmetry is a main
geometric difference that distinguishes some ``good'' $S$
from other ``bad'' $S$, it turns out that a large part of
the harmonic analysis on these groups can be worked out
regardless of this distinction.

A detailed account of this is given in our forthcoming
paper~\cite{7}.

We thank J.~Faraut and A.~Kor\'anyi for calling our
attention to the Lichnerowicz conjecture.

\heading 1. The extension $S$ of an H-type group\endheading

An H-type (or Heisenberg-type) algebra is a two-step
nilpotent Lie algebra $\germ n$ with an inner product
$\lan\ \ ,\ \ \ran$ such that if $\germ z$ is the center
of $\germ n$ and $\germ v =\germ z^\perp$, then the map
$J_Z:\germ v\rightarrow\germ v$ given by
$$
\lan J_ZX,Y\ran= \lan[X,Y].Z\ran
$$
for $X,Y\in \germ v$ and $Z\in \germ z$,
satisfies the identity $J_Z^2=-|Z|^2I$ for every
$Z\in\germ z$.

An H-type group is a connected and simply connected Lie
group $N$ whose Lie algebra is H-type.

Let $\germ s=\germ v \oplus \germ z \oplus \Bbb R T$ be
the extension of $\germ n$ obtained by adding the rule
$$
[T,X+Z]=\tfrac12 X+Z
$$
for $X\in \germ v$ and $Z\in \germ z$. Let $S=NA$ be the
corresponding connected and simply connected group
extension of $N$, where $A=\exp_S(\Bbb R T)$. We
parametrize the elements of $S$ in terms of triples
$(X,Z,a)\in \germ v \oplus \germ z \oplus \Bbb R^+$ by
setting
$$
(X,Z,a)= \exp_S(X+Z)\exp_S (\log a\,T)\ .
$$

The product on $S$ is then
$$
(X,Z,a)(X',Z',a')= (X+a^{1/2} X', Z+aZ' +\tfrac12
a^{1/2} [X,X'], aa')\ .
$$

We introduce the inner product
$$
\lan X+Z+tT|X'+Z'+t'T\ran = \lan X|X'\ran + \lan Z|Z'\ran
+tt'\
$$
on $\germ s$ and endow $S$ with the induced left-invariant 
Riemannian
metric.

\proclaim{Proposition 1 \cite{4,~5,~9}} Let $M=G/K$ be a
rank-one symmetric space of the noncompact type
different from $SO_e(n,1)/SO(n)$, and let $G=NAK$ be the
Iwasawa decomposition of $G$. Then $N$ is an H-type
group and the map $s\mapsto sK$ is an isometry from
$S=NA$ to $G/K$.
\endproclaim

From the classification of symmetric spaces one knows
that if the H-type group $N$ appears in the
Iwasawa decomposition of a rank-one symmetric space, then
the dimension of its center equals 1,~3, or~7. On the
other hand, there exist H-type groups with centers of
arbitrary dimensions~\cite{8}. Direct proofs can be
given~\cite{4,~5} to show that the space $S$ is symmetric
if and only if $N$ is an ``Iwasawa group.''

We can then conclude that {\it there are infinitely many
$S$ that are not symmetric}.

\heading 2. The volume element in radial normal
coordinates\endheading

Given unit elements $X_0\in \germ v$ and $Z_0\in \germ
z$, we denote by $S_{X_0,Z_0}$ the 4-dimensional subgroup
of $S$ generated by $X_0,\ Z_0,\ J_{Z_0}X_0$, and~$T$.
Clearly, $S_{X_0,Z_0}$ is also the extension of an H-type
group, precisely of the 3-dimensional Heisenberg group.

\proclaim {Proposition 2 \rm\cite{3}} $S_{X_0,Z_0}$ is a
totally geodesic submanifold of $S$, isometric to the
Hermitian symmetric space $SU(2,1)/S(U(2)\times U(1))$.
\endproclaim

This observation suggests two alternative realizations of
$S$~\cite{4,~7}:

\varroster
\item the ``Siegel domain'' model:
$$
D=\{(X,Z,t)\in \germ v \oplus \germ z \oplus \Bbb R :
t>\tfrac14|X|^2\}
$$
with the metric transported from $S$ via the map
$h(X,Z,a) =(X,Z, a+\frac14|X|^2)$;
\item the ``ball'' model\/\RM:
$$
B=\{(X,Z,t)\in \germ v \oplus \germ z \oplus \Bbb R :
|X|^2 +|Z|^2 +t^2<1\}
$$
with the metric transported from $D$ via the inverse of
the ``Cayley transform'' $C:B\rightarrow D$ given by
$$C(X,Z,t)= {1\over (1-t)^2+|Z|^2} \big(2(1-t+J_Z)X\ ,\ 
2Z\ ,\ 1-t^2-|Z|^2\big)\ .
$$
\endroster
The map $C^{-1}\circ h:S\rightarrow B$ maps the identity
element $e=(0,0,1)$ to the origin and its differential
at $e$ equals $\frac12 I$.

Let $m=\dim \germ v$, $k=\dim \germ z$. The number
$Q=\frac m2+k$ is called the homogeneous dimension of
$N$. It is easily checked that
$$
dm_L=a^{-Q-1}\,dX\,dZ\,da
$$
is a left Haar measure on $S$ and is also the Riemannian
volume element.

Applying Proposition~2, using the properties of the
Cayley transform in the symmetric Hermitian case, and
explicitely computing the Jacobians of $h$ and $C$, we
obtain the following result.

\proclaim {Theorem 1 \rm\cite{7}} 
{\rm(1)} In the metric on $B$ transported from $S$ via
$C^{-1}\circ h$, the geodesics through the origin are the
diameters and the distance $\rho$ from a point 
$(X,Z,t)\allowmathbreak \in B$
to the origin is a function of $r=
(|X|^2+|Z|^2+t^2)^{1/2}$, namely,
$$
\rho=\log\frac{1+r}{1-r}\ .
$$

{\rm(2)} Introducing polar coordinates $(r,\omega)$ on $B$
and denoting by $d\sigma(\omega)$ the surface measure on
$S^{m+k}$, the volume element on $B$ is given by
$$
2^{m+k+1}(1-r^2)^{-Q-1}r^{m+k}\,dr\,d\sigma(\omega) =
2^{m+k}\left(\cosh \frac\rho2\right)
^k\left(\sinh \frac\rho2\right)^{m+k}\,
d\rho\,d\sigma(\omega)\ .
$$
\endproclaim

The coordinates $(\rho,\omega)$ are the radial normal
coordinates on $B$ around the origin. Then (2) proves
that the volume density depends only on $\rho$. The
same is true for normal coordinates around any other
point, since $S$ and $B$ are isometric and $S$ is
obviously homogeneous. We can then conclude that

\proclaim{Corollary} For every H-type group $N$, $S=NA$
is a harmonic space.
\endproclaim

Consider now the Laplace-Beltrami operator $\Delta$ on
$S$. Its extension to the domain 
$$
\Cal D(\Delta) =\{f\in
L^2(S,dm_L):\Delta f\in L^2(S,dm_L)\}
$$ 
is selfadjoint
and positive. From the spectral resolution 
$$
\Delta=\int_0^{+\infty} \lambda\,dE_\lambda
$$
one constructs the heat semigroup
$$
T_t=e^{-t\Delta}=\int_0^{+\infty}
e^{-t\lambda}\,dE_\lambda\ .
$$

Since $\Delta$ is left-invariant,
$$
(T_tf)(x)=(f*p_t)(x)=\int_S f(y)p_t(y^{-1}x)\,dm_L(y)
$$
and $p(t,x)=p_t(x)$ is a smooth function on
$\Bbb R^+\times S$. In \cite{7} we also prove the
following strong-harmonicity property of $S$.
\proclaim{Theorem 2} The functions $p_t(x)$ depend only
on the geodesic distance from $x$~to~$e$.
\endproclaim

\heading 3. Dimensions of the nonsymmetric $S$\endheading

A complete list of H-type groups can be obtained from the
classification of Clifford modules~\cite{8}. For small
values of $k=\dim \germ z$ we write the dimensions of the
corresponding nonsymmetric $S$:
%
\roster
\item"" $k\qquad \dim S$
\item"" $1\qquad \text{---}$
\item"" $2\qquad 7+4n$
\item"" $3\qquad 12+4n$
\item"" $4\qquad 13+8n$
\item"" $5\qquad 14+8n$
\item"" $6\qquad 15+8n$
\item"" $7\qquad 24+8n$
\item"" $8\qquad 25+16n$
\endroster
where $n=0,1,\dots$.

\enddocument